\documentclass[12pt,a4paper]{article}
\usepackage[latin1]{inputenc}
\usepackage{amsmath}
\usepackage{amsfonts}
\usepackage{amssymb}
\usepackage{graphicx}
\usepackage{float}
\usepackage{authblk}
\usepackage{amsthm}
\usepackage{etoolbox}
\usepackage{caption}
\usepackage{tikz}
\usepackage{tikz-cd}
\usepackage{pst-node}
\usepackage{amsmath}
\usepackage{mathtools}
\usepackage{latexsym}
\usepackage{mathrsfs}
\usepackage{array} 
\usepackage{chngcntr}
\usetikzlibrary{arrows, matrix}
\numberwithin{equation}{section}
\usepackage[parfill]{parskip}
\newtheoremstyle{mystyle}
{7pt}
{5pt}
{\itshape}
{}
{\bfseries}
{}
{.5em}
{}

\theoremstyle{mystyle} 
\newtheorem{thm}{Theorem}[section]
\newtheorem{cor}[thm]{Corollary}
\newtheorem{lem}[thm]{Lemma}
\newtheorem{prop}[thm]{Proposition}
\theoremstyle{definition}
\newtheorem{defi}[thm]{Definition}
\newtheorem{exa}{Example}
\newtheorem{rem}{Remark}
\makeatletter
\begin{document}
\topmargin=-5mm
\textheight=220mm
\textwidth=190mm
\leftmargin=-4mm
\title{\Large \textbf{On Hilbert Functions of Points in Projective Space and Structure of Graded Modules}}
\author{Damas Karmel Mgani$^1$ and Makungu Mwanzalima $^2$\\
	$^1$Department of Mathematics, College of Natural and Applied Sciences, University of Dar es Salaam, TANZANIA. \\
	e-mail: d.mgani99@gmail.com\\
	$^2$Department of Mathematics, College of Natural and Applied Sciences, University of Dar es Salaam, TANZANIA.\\
	e-mail: fazakungu@gmail.com}

\date{}

 
\maketitle 
 
\begin{abstract}
\textsl{In this paper, we investigate the relationship between the Hilbert functions and the associated properties of the graded modules. To attain this, we construct the graded modules from the sets of points in projective space, $\mathbb{P}_k^n$ . We use a computer software package for algebraic computations Macaulay2 to study the Hilbert functions and the associated properties of the graded modules. Thereafter, we provide theoretical proofs of the results obtained from Macaulay2 and finally, we give illustrative examples to justify some of our results.}
\end{abstract}
\let\thefootnote\relax
\let\thefootnote\relax 

\section{Introduction and Preliminaries}
Let $R=k[x_0,...x_n]$ be a polynomial ring in $n+1$ independent variables over an algebraically closed field $k$, with $\text{deg}(x_i)=1$ for $i=1,...,n$ and let $X=\{P_1,...,P_s\}\subseteq \mathbb{P}_k^n$ be a finite set of points in a projective space with length $\textit{l}_{X}=s$, in this work we only consider sets of points with general position (random sets of points ). Now, we define an ideal $I_{X}$ of all functions vanishing on $X$, and $R_{X}=R/I_{X}$, the homogeneous coordinate ring of the points. Let $M$ be an $R_{X}$- module and denote by $H_i(M)$, for each $i\ge 0$, the dimension, as a vector space, of the graded component of degree $i$ of the graded module $M$. The function $i\mapsto H_i(M):=\dim_k M(i)=\binom{n+i}{i}$ is known as the Hilbert function of $M$ in degree $i$. The Hilbert function of a set of points in $\mathbb{P}_k^2$ is the basis for many questions about sets of points. To any set of points, we can associate an algebraic object which we call the coordinate ring. Then, the  Hilbert function is used to obtain, among other things, algebraic information about the coordinate ring and geometric information about the set of points \cite{Van}. 
The following definitions and preliminary results will be required to establish the results.  
\begin{thm}(Hilbert)\label{t1}
		If $M$ is a finitely generated graded module over $k[x_0,...x_n]$, then $H_i(M)$ agrees, for large $i$, with a polynomial of degree $\leq n$.
	\end{thm}
	The details of the proof of Theorem \ref{t1} can be obtained from \cite{Eisenbud},\cite{Li} and \cite{Atiyah}.
\begin{defi}\label{def12}
	Let $R_{X}$ be a Noetherian ring. A finite $R_{X}$- Module $M$ is torsionless if the map $\varphi_{M} : M\longrightarrow M^{**}$ is a monomorphism or injective, where $M^{**}=\text{Hom}_{R_{X}}(M^{*},R_{X})$ is a double dual of $M=\text{Hom}_{R_{X}}(J,R_{X})$. The use of the canonical homomorphism of a module into its double dual was initiated by Bass \cite{Bass}.
\end{defi}
\begin{defi}\label{def2}
	An $R_{X}$- module $M$ is torsionless if and only if for each $0\neq m\in M$ there exists $f\in M^{*}$ such that $f(m)\neq 0$.
\end{defi}
\begin{defi}
	If  	
$	0{\longrightarrow}F_0\stackrel{f}{\longrightarrow}F_{1} \stackrel{g}{\longrightarrow}F_2\longrightarrow 0$
  is an exact sequence of $R_{X}$- modules and $R_{X}$- homomorphisms, so is  
$	0{\longrightarrow}F_2^{*}\stackrel{g^{*}}{\longrightarrow}F_1^{*} \stackrel{f^{*}}{\longrightarrow}F_0^{*}. $
\end{defi}
\begin{defi}\label{def4} An ideal $I$ in $R$ is said to be  regular if and only if it consists entirely of regular elements (or non-zero divisors) of $R$.
\end{defi}
\begin{rem}
	\begin{enumerate}
		\item[(A)]  It is easy to see that a regular ideal $I$ in $R$ is itself a regular ring. For if $u\in I$, $\exists x \in R$ such that $uxu=u$. Next, $uxuxu=u$ and $xux\in I$. Hence, $u$ is regular in the ring $I$.
		\item[(B)] Since our rings are commutative, then $uxu=u(ux)=u^2x=u$.
	\end{enumerate}
	
\end{rem}
\begin{prop} \cite{Mwanzalima}\label{prop1} Let $J\subseteq R_{X}$ be a non-zero homogeneous ideal of $R_{X}$ and $J\hookrightarrow R_{X}$ an inclusion map, then $\text{Hom}_{R_{X}}(J,R_{X})$ is an $R_{X}$- module of homomorphisms.
\end{prop}
\begin{lem} \emph{\cite{Mwanzalima}}\label{lem3} Let $J\subseteq R_{X}$ be any homogeneous ideal. Then, the map $\mu : J\longrightarrow \text{Hom}_{R_{X}}(\text{Hom}_{R_{X}}(J,R_{X}),R_{X})$ is injective.
\end{lem}
\begin{prop}\cite{Mwanzalima}\label{prop2} If $J\subseteq R_{X}$ is a homogeneous ideal and satisfies the conditions, $\dim_k \text{Hom}_{R_{X}}(J,R_{X})_i=\dim_k (R_{X})_i=s$ for $i\geq\delta_X$, then $\text{Hom}_{R_{X}}(J,R_{X})_i$ as an $R_{X}$- module of homomorphisms is isomorphic to some homogeneous ideal of $R_{X}$. $\delta_{X}$ is the lowest degree for which the Hilbert function of $X$ equals the number of points $s$.
\end{prop}
\newpage
\section{Main Results}
In this section we give proofs of our main results.
\begin{lem} \label{lem2}If $J\subseteq R_{X}$ is a regular homogeneous ideal, then the map $\theta_g: \text{Hom}_{R_{X}}(J,R_{X})\longrightarrow R_{X}$ is injective if and only if $g$ is a non-zero divisor on $R_{X}$.
	\begin{proof}
		$(\Longleftarrow)$: Let $g\in J$ be a non-zero divisor on $R_{X}$. Then for all $f\in J$ we have
		\begin{align}\phi(fg) =g\phi(f)=f\phi(g).    \nonumber
		\end{align} 
		If $\phi(g)=0$ then $g\phi(f)=0$. So, we have $\phi(f)=0$ for all $f\in J$ giving $\phi=0$. This means that $\theta_g(\phi)=0$, implies that $\phi =0$. Hence $\theta_g$ is injective.\\
		$(\Longrightarrow):$ Suppose $g\in J$ is a zero divisor. Then there exists $h\neq 0$ such that $gh=0$. Let $\phi \in \text{Hom}_{R_{X}}(J,R_{X})$ be defined by $\phi(f)=fh$ for all $f\in J$. Then we have $\phi(g)=gh=0$. That is, $\theta_g(\phi)=\phi(g)=0$
		and we get that $\phi \in ker(\theta_g)$. But $\phi \neq 0$. Take $f\in J$ to be a non-zero divisor on $R_{X}$. Then $\phi(f)=fh\neq 0$, because $h\neq$ 0. Thus $\theta _g$ is not injective.
	\end{proof}
\end{lem}
\begin{prop} If $h\in J$ is a non-zero divisor, then there exists at least one homomorphism $\text{Hom}_{R_{X}}(J,R_{X})\longrightarrow R_{X}$ such that $R_{X}/J$ is Artinian.
	\begin{proof} Since $J$ contains a non-zero divisor $h$, it follows from Definition \ref{def4} that $J$ is a regular homogeneous ideal of $R_{X}$, and from Lemma \ref{lem2},\\ $\text{Hom}_{R_{X}}(J,R_{X})\longrightarrow R_{X}$ exists and is injective. Next, we  know that $R_{X}=R/I_{X}$, suppose $R$ is an Artinian ring with ideal $I_{X}$, then $R/I_{X}$ is also Artinian (a quotient of an Artinian ring is Artinian). Now, it is easy to see that $R_{X}/J$ is Artinian. 
		Alternatively, since  $R_{X}/J$ is of dimension zero ($\dim (R_{X}/J) =0$) and the Hilbert function $H(R_{X}/J, i)=0$ for $i$ large. Hence, $R_{X}/J$ is Artinian.	
	\end{proof}
		\begin{rem}
			The homomorphism is of the same degree as that of the minimal generator(s) of $J$.
		\end{rem}
\end{prop}
	\begin{cor}\label{cor1} 
		Let $X \subseteq \mathbb{P}_k^2$ be a set of points  with  $R_{X}/J$ Artinian. Then the ideal $J$ containing at least one generator of $I_{X}$ is generated by homogeneous elements of initial degree $i\leq \Omega$, where $\Omega$ is a Socle degree of $R_{X}/J$.
		\begin{proof}
			Let 
			\begin{eqnarray}
			0\longrightarrow \bigoplus\limits_{i=1}^{P}R(-d_i)\longrightarrow \bigoplus\limits_{i=1}^{p+1}R(-c_i)\longrightarrow J\longrightarrow 0.		 
			\end{eqnarray}
			be a minimal free resolution of $J$, where $c_1\leq c_2\leq ... \leq c_{p+1}$ are the degrees of the minimal set of generators of $J$ and $d_1\leq d_2\leq...\leq d_p$ the degrees of a minimal set of generators for the syzygy module.  For $i=1,...,p+1$ we write $\lambda_j=|\{c_i|c_i=j\}|$
			Let $\lambda _j$ be the number of the integers $c_i$ equal to $j$. Let $H_i$ be the Hilbert function of $R_{X}/J$, that is, $H_{i}=Hilb(R_{X}/J,i)$ and $q=min\{i\in \mathbb{N}| H_{i}=0\}$. Since $H_i(\Omega)\neq 0$ and $H_{j}=0$ for all $j\geq \Omega$, the highest value of $j$ for which $\lambda_j$ could be non zero is $\Omega$. 	
		\end{proof}
	\end{cor}
	In the next proposition (Proposition \ref{tor1}), we prove that any submodule of a torsionless $R_{X}$- module is torsionless over Noetherian rings.
	\begin{prop}\label{tor1}
		Let $R_{X}$ be a Noetherian ring and $M$ a torsionless $R_{X}$- module. If $N\subset M$ is an $R_{X}$- submodule, then $N$ is torsionless.
	\end{prop}
	\begin{proof}
		Since $N\subset M$, there is an injection 
		$0\longrightarrow N \longrightarrow M.$\\
		From definition \ref{def12} , we have an injective map $\varphi_{M} : M\longrightarrow M^{**}$.
		Now, we construct $\varphi_{N} : N\longrightarrow N^{**}$, where $N^{**}$ is the double dual of $N$.
		Then, we get the following commutative diagram
		\[
		\begin{tikzcd}
		0 \arrow[r," "]
		& N \arrow[r,"f"] \arrow[d, ,"\varphi_{N}" ] & M\arrow[d, ,"\varphi_{M}" ]\\
		& N^{**} \arrow[r,"f**"]& M^{**}
		\end{tikzcd}
		\]\\
		where $f^{**}=(f^*)^{*}$. \\
		From the commutative diagram one has that,
		$\varphi_{M}\circ f$ is injective (composition of injective maps is injective), and $f^{**}\circ \varphi_{N}$ is also injective.
		It follows  that $\varphi_{N}$ is injective. Hence, $N$ as an $R_{X}$-  submodule of $M$ is torsionless.		
	\end{proof}

\begin{exa}

		All the computations in this section were done using the computer algebra package Macaulay2. We describe  the relationship between the initial degree of the miniminal generator(s) of a homogeneous ideal $J$ and the Socle degree of an Artinian algebra $R_{X}/J$. We investigate through different sets of points in the projective plane $\mathbb{P}_k^2$. Let  $X_{1}=\{P_1, P_2\}, X_{2}=\{P_1, P_2, P_3\}$ and $X_{3}=\{P_1, P_2, P_3, P_4\}$ be sets with $2, 3$ and $4$ random points respectively (we call this, the first group of points).
		We construct the homogeneous ideals $I_{X_{1}}, I_{X_{2}}$  and $I_{X_{3}}$  and compute the Betti diagrams as shown hereunder:
		\begin{table}[H]			
			\centering 
			\caption{Betti diagram for $I_{X_{1}}$} 
			\label{tab:table pt2}
			\begin{tabular}{c| c c c c c}
				& 1& 2 & 1 & \\ [0.5ex]
				
				\hline
				0: & 1& 1&-&\\
				1: &-& 1& 1 & \\   						
			\end{tabular}			
		\end{table}  
		$I_{X_{1}}$ is generated by $1$ form of degree $1$ and $1$ form of degree $2$.
		\begin{table}[H]			
			\centering 
			\caption{Betti diagram for $I_{X_{2}}$} 
			\label{tab:table pt3}
			
			\begin{tabular}{c| c c c c c}
				& 1& 3 & 2 & \\ [0.5ex]
				
				\hline
				0: & 1& -&-&\\
				1: &-& 3& 2 & \\   						
			\end{tabular}			
		\end{table} 
		$I_{X_{2}}$ is generated by $3$ forms of degree $2$.
		\begin{table}[H]		
			\centering 
			\caption{Betti diagram for $I_{X_{3}}$} 
			\label{tab:table pt33}
			\begin{tabular}{c| c c c c c}
				& 1& 2 & 1 & \\ [0.5ex]
				\hline
				0: & 1& -&-&\\
				1: &-& 2& - & \\  
				2: & - & -& 1 		
				
			\end{tabular}			
		\end{table}  
		$I_{X_{3}}$ is generated by $2$ forms of degree $2$.\\
		We constuct $J_{X_{1}}, J_{X_{2}}$ and $J_{X_{3}}$ as homogeneous ideals containing at least a degree of the generator(s) of $I_{X_{1}}, I_{X_{2}}$ and $I_{X_{2}}$ respectively.
		Taking $J_{X_{1}}=(x_0+x_1+x_2, x_0^2+x_1^2+x_2^2-x_0x_1+x_1x_2)$, $J_{X_{2}}=(x_0+x_1+x_2, x_0^2+x_1^2+x_2^2+x_0x_1+x_1x_2-x_0x_2)$,  and $J_{X_{3}}=(x_0+x_1+x_2, x_0^2+x_1^2+x_2^2-x_0x_1+x_1x_2)$. $Hilb(R_{X}/J_{X_{1}})=\{1,1,0,0,0,...\}$, $Hilb(R_{X}/J_{X_{2}})=\{1,2,0,0,0,...\}$, and $Hilb(R_{X}/J_{X_{3}})=\{1,2,0,0,0,...\}$. Since all the Hilbert functions are zero for higher degrees, thus $R_{X}/J_{X_{1}}, R_{X}/J_{X_{2}}$, and $R_{X}/J_{X_{3}}$ are Artinian, both of socle degree $1$. We see that, for sets of points $X_1, X_2$ and $X_3$, the initial degree of the generator(s) of  the ideals $J_{X_{1}}, J_{X_{2}}$ and $J_{X_{3}}$ is the same as the socle degree of the corresponding coordinate ring.\\
		For $X=\{P_1,...,P_s\}$, where $s\geq 4$, situation is different, the socle degree increases as does the initial degree of $J$. Again, this will be illustrated below for some sets of points, and later on we give the general  behavior for any finite set of points in the projective plane $\mathbb{P}_k^2$.\\
		Let $X_{6}=\{P_1,...,P_7\}$ be a set with $7$ random points, $I_{X_{6}}$ is generated by $3$ form of degree $3$ and the corresponding Betti diagram is as shown below,
		\begin{table}[H]		
			\centering 
			\caption{Betti diagram for $I_{X_{6}}$} 
			\label{tab:table pt4}
			
			\begin{tabular}{c| c c c c c}
				& 1& 3 & 2 & \\ [0.5ex]
				
				\hline
				0: & 1& -&-&\\
				1: &-& -& - & \\  
				2: & - & 3& 1 & \\
				3: &- & -& 1		
				
			\end{tabular}			
		\end{table}   
		Let $J_{X_{6}}=(x_0^2+x_1^2+x_2^2+x_0x_1+x_1x_2, x_0^3+x_1^3+x_2^3-x_0^2x_1+x_1x_2^2+x_0x_1x_2)$, we have $Hilb(R_{X}/J_{X_{6}})=\{1,3,5,3,0,0,0,...\}$. $J_{X_{6}}$ is of initial degree $2$ and the highest degree for which its Hilbert function is non-zero is $3$. Thus, $R_{X}/J_{X_{6}}$ is Artinian of Socle degree $3$. Equivalently, $
		Socle \enskip degree = \text{deg}(\text{min}\{G_{J_{X_{6}}}\})+1
		=2+1=3$.
		This holds for $5\leq s \leq 9$, ( the second group ).
		We also observe that, for the third group, $10\leq s \leq 16$, $Socle \enskip degree = \text{deg}(\text{min}\{G_{J_{X}}\})+2$.
		For the fourth group, $17\leq s \leq 25$, $	Socle \enskip degree = \text{deg}(\text{min}\{G_{J_{X}}\})+3$
		and so on.
		%
		%
		%
		Where $G_{J_{X}}$, is the set of generator(s) of $J_{X}$.\\ Generally, let $ \xi $ denote the collection of distinct sets of points with the same socle degree. Then, the number of sets in this collection is obtained by the following formula 
		 \begin{eqnarray}\nonumber
		\xi = 2n + 1, \quad \forall n\geq 1.
		\end{eqnarray}
	\end{exa}
In the following result we wish to investigate what happens to quotient rings $R_{X}/I$, $R_{X}/I^{*}$ and $R_{X}/(I+I^*)$ , when $I$ and $I^{*}$ are monomial ideals. 
An ideal $I\subseteq k[x_{0},...,x_{n}]$ is termed to be monomial  if it is generated by monomials of the form $x^a=x_1^{a_0}... x_n^{a_n}, a=(a_0,..., a_n)$. To find a monomial ideal $I$ that contains at least one monomial from each polynomial of $I_X$ it is important to look for a method  of choosing one monomial from each polynomial of $R$. For our case, we choose leading terms  from each polynomial in $I_X$, and use the chosen monomials or terms to generate $I$. We use $I$ and $I^{*}$ to be the monomial ideals generated from $I_{X_1}$ and $I_{X_2}$ respectively. 
\begin{lem}
	Let $X_1,X_2\subseteq \mathbb{P}_k^2$ be two sets of random points, $I_{X_1}$ and $  I_{X_2}$ be homogeneous ideals generated by all forms vanishing on $X_1$ and $X_2$ respectively. Let $I$ be a monomial ideal of $I_{X_1}$  and, $I^*$ a monomial ideal of $I_{X_2}$. Then, the quotients rings $R_{X}/I, R_{X}/I^*$ and $R_{X}/(I+I^*)$ are Artinian.
	\begin{proof}
		In order to prove that $R_{X}/I$ is Artinian, we need to show that for any two distinct prime ideals $I$ and  $I^\prime$ there does not exist an ideal properly containing the other, that is $I\subset I^\prime$ or $I^\prime\subset I$. We begin by assuming that $I$ is a prime ideal of $R_{X}$, then we know that $R_{X}/I$ is an integral domain. For $R_{X}/I$ to be Artinian we must have the descending chain condition (d.c.c). Suppose  $x\in R_{X}/I$, $x$ is a nonzero element, by d.c.c it follows that $(x^n)=(x^{n+1})$ for some $n$. Now, for some $r\in R_{X}/I$ we can have $x^n=x^{n+1}r$. Now, the cancellation laws hold to $x^n$ since $R_{X}/I$ is an integral domain and $x$ is nonzero, thus $xr=1$. Therefore, $r$ is a multiplicative inverse of $x$ and $R_{X}/I$ is a field, hence $I$ is a maximal ideal. So there exists no other ideal $I^\prime$ of $R_{X}$ such that $I\subset I^\prime$. Rings having this property are refered to as rings of krull dimension $0$. This shows that the quotient ring $R_{X}/I$ is Artinian, so is $R_{X}/I^*$.\\
		To prove that $R_{X}/(I+I^*)$ is Artinian. We know that 
		\begin{align}\nonumber
			I+I^*&=<\text{Monomials in both } I_{X_1} \text{and} \, I_{X_2}>\nonumber\\
			&=I_{X_1}\cap I_{X_2} \nonumber\\
			&=I_{X_1\cup X_2} \label{dam}
		\end{align}
		Where $I_{X_1\cup X_2}$ is the monomial ideal of the union of two sets of points $X_1$ and $X_2$.
		Now, it follows that 
		\begin{align}\nonumber
			R_{X}/(I+I^*)& \cong R_{X}/(I_{X_1}\cap I_{X_2})=R_{X}/(I_{X_1\cup  X_2})\nonumber
		\end{align} 
		But either $I_{X_1}\cap I_{X_2} \subseteq I_{X_1}$ or $ I_{X_1}\cap I_{X_2} \subseteq I_{X_2}$. We get that
		\begin{align}\nonumber 
			R_{X}/(I_{X_1}\cap I_{X_2}) &\subseteq R_{X}/I_{X_1} \cup R_{X}/I_{X_2} \nonumber 
		\end{align}
		Since $R_{X}/I_{X_1}$ and  $R_{X}/I_{X_2}$ are both Artinian,  implies that $R_{X}/(I_{X_1}\cap I_{X_2})$ is also Artinian. Hence $R_{X}/(I+I^*)$ is Artinian.
	\end{proof}
\end{lem}	
We illustrate the discussion above with the following example:
\begin{exa}

	Let $X_1=\{P_1,...,P_{15}\}$ be a set with $15$  points and $I_{X_1}$ the ideal vanishing on $X_1$. The associated monomial ideal of $I_{X_1}$ is $I=(x_1^5,x_0x_1^4,x_0^2x_1^3,\\x_0^3x_1^2,x_0^4x_1,x_0^5)$. The Hilbert function  of $R_{X}/I$, $\text{Hilb}(R_{X}/I)=$ $\{1,3,6,10,15,9,\\2,0,0,0,...\}$. Take, $X_2=\{P_1,...,P_{21}\}$ and $I_{X_2}$, an ideal vanishing on $X_2$ with corresponding monomial ideal $I^*=(x_0^6,x_0^5x_1,x_0^4x_1^2,x_0^3x_1^3,x_0^2x_1^4,x_0x_1^5,x_1^6)$,  $\text{Hilb}(R_{X}/I^*)=\{1,3,6,10,15,15,8,0,0,0,...\}$. Then, using Macaulay 2 one can verify that we have $I+I^*$ as an ideal of $R_{X}$ with Hilbert function, $\text{Hilb}(R_{X}/(I+I^*))=\{1,3,6,15,9,2,0,0,...\}$. Since all the Hilbert functions become zero for large degree $i$, then the quotient rings, $R_{X}/I, R_{X}/I^*$ and $R_{X}/(I+I^*)$ are all Artinian.
\end{exa}

\end{document}